\documentclass{article}
\usepackage{amssymb,amsmath,bm,hyperref}
\usepackage{amssymb,amsmath,bm,geometry,graphics,graphicx,url,color,inputenc}
\def\no{\noindent}
\def\pmatrix{\left(\begin{array}}
\def\endpmatrix{\end{array}\right)}

\def\RR{\mathbb{R}}

\def\B{{\cal B}}

\def\F{{\cal F}}

\def\I{{\cal I}}

\def\P{{\cal P}}

\def\dd{\mathrm{d}}

\def\diag{\mathrm{diag}}

\newtheorem{theo}{Theorem}
\newtheorem{lem}{Lemma}

\newtheorem{rem}{Remark}
\newtheorem{defi}{Definition}
\def\proof{\noindent\underline{Proof}\quad}
\def\QED{\mbox{~$\Box{~}$}}

\def\bfc{{\bm{c}}}

\def\bfe{{\bm{e}}}

\def\bfzero{{\bm{0}}}

\def\bfphi{{\bm{\phi}}}

\def\om{\omega}

\begin{document}

\title{Arbitrary high-order methods for one-sided\\ direct event location in discontinuous differential problems\\ with nonlinear event function}

\author{Pierluigi Amodio\footnote{Dipartimento di Matematica, Universit\`a di Bari, Italy. \qquad\, \url{{pierluigi.amodio,felice.iavernaro}@uniba.it}} \and Luigi Brugnano\footnote{Dipartimento di Matematica e Informatica ``U.\,Dini'', Universit\`a di Firenze, Italy. \quad\url{luigi.brugnano@unifi.it}}  \and  Felice Iavernaro$^*$}


\maketitle

\begin{abstract}  In this paper we are concerned with numerical methods for the one-sided event location in discontinuous differential problems, whose event function is nonlinear (in particular, of polynomial type). The original problem is transformed into an equivalent Poisson problem, which is effectively solved by suitably adapting a recently devised class of energy-conserving methods for Poisson systems. The actual implementation of the methods is fully discussed, with a particular emphasis to the problem at hand. Some numerical tests are reported, to assess the theoretical findings.

\medskip
\no{\bf Keywords:~} discontinuous ODEs, Poisson problems, Line Integral Methods, Hamiltonian Boundary Value Methods, HBVMs, PHBVMs, EPHBVMs.

\medskip
\no{\bf MSC:~} 65L05, 65P10.

\end{abstract}

\section{Introduction}\label{intro} 
In some applications, one faces the problem of solving {\em discontinuous ODE problems}, namely, problems in the form:
$$\frac{\dd}{\dd\tau} x = \left\{ \begin{array}{cc} f(x), & \mbox{~if~} g(x)\le0,\\[2mm] \phi(x), &\mbox{otherwise,}\end{array}\right. \qquad x(0)=x_0\in\RR^n,$$ 
with ~$g:\RR^n\rightarrow \RR$~ a suitably regular function, called {\em event function}, dividing $\RR^n$ into two regions where the vector field is defined in different ways. The vector field need not be continuous on the boundary set
\begin{equation}\label{Sset}
\Sigma = \left\{ x\in\RR^n\,:\,g(x)=0\right\}.
\end{equation} 
Hereafter, we shall refer to the set $\Sigma$ as to the {\em event set} and, moreover, we shall assume, without loss of generality,  that $g(x_0)<0$.  
This problem has been studied in \cite{DL2015} (see also \cite{DL2012,LM2017} and references therein), in the case where $g(x)$ is linear or, at most, quadratic: here, the authors define a {\em direct} method for finding the {\em event point}, namely, the first point $x^*$ of the trajectory belonging to $\Sigma$. We refer to the references in \cite{DL2015} for relevant applications where solving such a problem is needed. 

In this paper, we consider the more general case where $g$ is a general polynomial. Consequently, given the ODE-IVP
\begin{equation}\label{lope_ode}
\frac{\dd}{\dd \tau} x = f(x), \qquad x(0) = x_0\in\RR^n,
\end{equation}
 the problem at hand is that of determining the point $x^*$, on the solution trajectory, such that 
\begin{equation}\label{lope_h}
g(x^*)=0\in\RR, 
\end{equation}
where $g$ is a polynomial (however, we shall also sketch  the case of a general, sufficiently smooth, function). Hereafter, it is assumed that:
\begin{itemize}
\item $f$ is sufficiently smooth;

\item at the considered initial point,
\begin{equation}\label{lope_H}
g(x_0) = -\bar H<0,
\end{equation}

\item along the solution of (\ref{lope_ode}),
\begin{equation}\label{lope_delta}
\frac{\dd}{\dd \tau}g(x) = \nabla g(x)^\top f(x)\ge \delta>0,
\end{equation}
which implies that the set $\Sigma$ in (\ref{Sset}) is attractive for the trajectory starting at $x_0$. More precisely, any trajectory starting at a point $x_0$ satisfying (\ref{lope_H}) and (\ref{lope_delta}), will reach $\Sigma$ in a {\em finite time} $\tau^*$ (depending on the starting point).\footnote{In fact, by virtue of (\ref{lope_delta}), one obtains that $\tau^*\le\bar H/\delta$.}
\end{itemize}
As stated above, we shall refer to the vector $x^*$ satisfying (\ref{lope_h}) as the {\em event point}. Its detection is significant in many applications, where it is mandatory that the set $\Sigma$ is not crossed, but just reached by the trajectory \cite{DL2015}. Therefore, it makes sense to impose the same requirement to a numerical method of approximation. 

With these premises, the structure of the paper is as follows: in Section~\ref{Pform} we cast the problem (\ref{lope_ode})--(\ref{lope_H}) in Poisson form; in Section~\ref{EPHBVM} we recall the basic facts about the solution procedure for the new formulation, which is duly adapted for the problem at hand; in Section~\ref{numtest} we provide some numerical tests; at last, in Section~\ref{fine} a few conclusions are given.

\section{Poisson formulation}\label{Pform}

As previously observed, the numerical procedures studied in \cite{DL2015,LM2017} allow to effectively solve the problem (\ref{lope_ode})--(\ref{lope_H}) when $g$ is a {\em linear} or a {\em quadratic} function, by using a suitable reparametrization of time. In particular, we shall use a reparametrization akin to that used in \cite{DL2015}, i.e.,
\begin{equation}\label{lope_s}
\om = g(x(\tau)) +\bar H,
\end{equation}
which allows solving the problem in the interval $\omega\in[0,\bar H]$, due to (\ref{lope_H}) and to the monotonicity of $g$ along the solution of (\ref{lope_ode}) (see (\ref{lope_delta})). Consequently, by introducing the augmented state vector
\begin{equation}\label{lope_y}
y = \pmatrix{c} x\\ \om\endpmatrix\in \RR^m, \qquad m:=n+1,
\end{equation}
and using the independent variable 
\begin{equation}\label{teqs}
t=\om, 
\end{equation}
we obtain the augmented system, equivalent to (\ref{lope_ode}),
\begin{equation}\label{lope_odey}
\dot y = G(y), \qquad t\in[0,\bar H], \qquad y(0) = y_0 := \pmatrix{c} x_0\\ 0\endpmatrix,
\end{equation}
where, hereafter, $\dot y$ will denote the derivative w.r.t. $t$, and (see (\ref{lope_ode}) and (\ref{lope_y}))
\begin{equation}\label{lope_G}
G(y) := \pmatrix{c} \frac{f(x)}{\nabla g(x)^\top f(x)}\\[1mm] 1\endpmatrix.
\end{equation}
Clearly, because of (\ref{lope_s}), the problem (\ref{lope_odey}) has the scalar invariant
\begin{equation}\label{lope_I}
H(y) := g(x)-\om +\bar H.
\end{equation}
In fact, one has (see (\ref{lope_y}) and (\ref{lope_G})):
\begin{eqnarray}\nonumber
\dot H(y) &=& \nabla H(y)^\top\dot y ~=~  \pmatrix{cc} \nabla g(x)^\top & -1\endpmatrix G(y) \\[1mm] \label{derH}
&=& \pmatrix{cc} \nabla g(x)^\top & -1\endpmatrix
\pmatrix{c} \frac{f(x)}{\nabla g(x)^\top f(x)}\\[1mm] 1\endpmatrix ~=~ 1-1 ~=~ 0.
\end{eqnarray}
Since (see (\ref{lope_H}))
$$H(y_0) = g(x_0)-\om(0) +\bar H= -\bar H+\bar H=0,$$
at $t=\om=\bar H$ it will be 
\begin{equation}\label{ybH}
y(\bar H) = \pmatrix{c} x^*\\ \bar H\endpmatrix,\qquad x^*:=x(\bar H),
\end{equation}
such that
$$0~=~H(y(\bar H)) ~=~ g(x(\bar H)) - \omega(\bar H) + \bar H ~\equiv~  g(x^*)-\bar H + \bar H ~=~ g(x^*),$$ 
thus recovering the event point $x^*\in\Sigma$.

The novelty of the present paper is that of deriving procedures able to reach the event point in a finite number of steps, in the case where $h\in\Pi_\nu$.\footnote{As is usual, $\Pi_\nu$ denotes the vector space of polynomials of degree not larger than $\nu$.} The basic idea is that of transforming the original system (\ref{lope_odey}) into an equivalent Poisson problem:
\begin{equation}\label{poisson}
\dot y = B(y)\nabla H(y), \quad t\in[0,\bar H],\qquad y(0)=y_0\in\RR^m, \qquad B(y)^\top = -B(y).
\end{equation}
The following result, based on \cite{QC1996}, holds true.

\begin{theo}\label{equivlenti}
Problems (\ref{lope_odey})--(\ref{lope_G}) and (\ref{poisson}) are equivalent, and possess the invariant $H(y)\equiv0$, provided that the skew-symmetric matrix $B(y)$ is defined as follows:
\begin{equation}\label{By}
B(y) = \frac{G(y)\nabla H(y)^\top - \nabla H(y)G(y)^\top}{\|\nabla H(y)\|_2^2}. 
\end{equation}
\end{theo}
\proof 
First of all, we observe that matrix (\ref{By}) is well defined, since $\|\nabla H(y)\|_2^2>1$ (see (\ref{lope_delta}),  (\ref{lope_y}), and (\ref{lope_I})). Next, for the problem (\ref{poisson}) one has
$$\dot H(y) = \nabla H(y)^\top\dot y =  \nabla H(y)^\top B(y)\nabla H(y) = 0,$$
due to the fact that $B(y)$ is skew-symmetric. Moreover, since $\nabla H(y)^\top G(y)=0$ (see (\ref{lope_G}) and (\ref{derH})), one has:
\begin{eqnarray*}
B(y) \nabla H(y) &=&  \frac{G(y)\nabla H(y)^\top - \nabla H(y)G(y)^\top}{\|\nabla H(y)\|_2^2}\nabla H(y)\\
&=&  \frac{G(y)\overbrace{\nabla H(y)^\top\nabla H(y)}^{=\|\nabla H(y)\|_2^2} - \nabla H(y)\overbrace{G(y)^\top\nabla H(y)}^{=0}}{\|\nabla H(y)\|_2^2}\\[2mm]&=&G(y).~\QED
\end{eqnarray*}

When matrix $B(y)$ is constant, as in the case of Hamiltonian problems,
\begin{equation}\label{Ham}
 \dot y = J\nabla H(y), \qquad t>0, \qquad y(0)=y_0, \qquad J^\top=-J,
 \end{equation}
then $H(y)$ is referred to as the {\em energy}, and its conservation can be effectively and efficiently obtained by solving problem (\ref{Ham}) via {\em Hamiltonian Boundary value Methods (HBVMs)}, a class of energy-conserving Runge-Kutta methods for Hamiltonian problems (see, e.g., \cite{BIT2009,BIT2010,BIT2011,BIT2012-1,BIT2012,BFCI2014,BIT2015} and the monograph \cite{LIMbook2016}, see also the review paper \cite{BI2018}). Nevertheless, in the case where the problem is not Hamiltonian, HBVMs are no more energy-conserving. When $B(y)=-B(y)^\top$ is not constant, problem (\ref{poisson}) is a particular instance of a Poisson problem. 
This motivates the present paper, where a recently-derived energy-conserving variant of HBVMs for Poisson problems \cite{ABI2022} will be suitably adapted for solving problem (\ref{poisson})-(\ref{By}).

For sake of completeness, we mention that the numerical solution of Poisson problems has been tackled by following many different approaches (see, e.g., \cite[Chapter\,VII]{GNI2006} and references therein). More recently, it has been considered in \cite{CH2011}, where an extension of the AVF method is proposed, and in \cite{BCMR2012,BGI2018}, where a line integral approach has been used instead. Functionally fitted methods have been proposed in \cite{M2015,MHW2022,WW2018}.

\section{Poisson HBVMs and their enhanced version}\label{EPHBVM}
Let us sketch the {\em Poisson HBVMs (PHBVMs)} methods defined in \cite{ABI2022}, which will be later slightly modified for the problem at hand.
Since we deal with one-step methods, we can consider the solution of problem (\ref{poisson}) on the interval $[0,h]$, with $h>0$ the timestep.  The basic idea is that of expanding the vector field (\ref{poisson}) along the orthonormal Legendre polynomial basis,
\begin{equation}\label{orto}
P_i\in\Pi_i, \qquad \int_0^1 P_i(\xi)P_j(\xi)\dd \xi = \delta_{ij}, \qquad i,j=0,1,\dots,
\end{equation}
with $\delta_{ij}$ the Kronecker symbol. In so doing, with similar steps as in \cite{ABI2022}, by considering the expansions
\begin{eqnarray}\nonumber
\nabla H(y(ch)) = \sum_{j\ge0} P_j(c) \gamma_j(y), &&P_j(c)B(y(ch))  = \sum_{i\ge0} P_i(c) \rho_{ij}(y), \qquad~  c\in[0,1],\\[-3mm] \label{gammaro}
\\[-3mm]\nonumber
\gamma_j(y)     = \int_0^1 P_j(\xi)\nabla H(y(\xi h))\dd\xi, &&
\rho_{ij}(y)         = \int_0^1 P_i(\xi)P_j(\xi)B(y(\xi h))\dd\xi, \quad~ i,j=0,1,\dots,
\end{eqnarray}
one obtains:
\begin{eqnarray}\label{infiy1}
\dot y(ch) &=& B(y(ch))\nabla H(y(ch)) = B(y(ch))\sum_{j\ge0} P_j(c)  \gamma_j(y)\\ \nonumber
&=& \sum_{j\ge0} P_j(c) B(y(ch)) \gamma_j(y) =  \sum_{i,j\ge 0} P_i(c)\rho_{ij}(y) \gamma_j(y), \qquad c\in[0,1],
\end{eqnarray}
from which one derives that the solution of (\ref{poisson}) can be formally written as:
\begin{equation}\label{infiy}
y(ch) = y_0 + h\sum_{i,j\ge0} \int_0^cP_i(\xi)\dd \xi\rho_{ij}(y) \gamma_j(y), \qquad c\in[0,1].
\end{equation}
In particular, by considering (\ref{orto}) and that $P_0(\xi)\equiv 1$, from which $\int_0^1P_i(\xi)\dd\xi=\delta_{i0}$ follows, one has:
\begin{eqnarray}\nonumber 
y(h) &=& y_0 + h\sum_{j\ge 0}\rho_{0j}(y) \gamma_j(y) \\ \label{infiyh}
&\equiv& y_0 + h\sum_{j\ge0} \int_0^1 P_j(\xi)B(y(\xi h))\dd \xi\int_0^1 P_j(\xi) \nabla H(y(\xi h))\dd \xi. 
\end{eqnarray}
In order to obtain a polynomial approximation of degree $s$ to $y$, it suffices to truncate the two infinite series in (\ref{infiy1}) after $s$ terms:
\begin{equation}\label{sigma1}
\dot \sigma(ch) =   \sum_{i,j=0}^{s-1} P_i(c)\rho_{ij}(\sigma) \gamma_j(\sigma), \qquad c\in[0,1],
\end{equation}
with $\rho_{ij}(\sigma)$ and $\gamma_j(\sigma)$ defined according to (\ref{gammaro}) by formally replacing $y$ by $\sigma$.
Consequently, (\ref{infiy}) becomes
\begin{equation}\label{sigma}
\sigma(ch) = y_0 + h\sum_{i,j=0}^{s-1} \int_0^cP_i(\xi)\dd \xi\rho_{ij}(\sigma) \gamma_j(\sigma), \qquad c\in[0,1],
\end{equation}
providing the approximation
\begin{eqnarray}\nonumber
y_1&:=&\sigma(h) = y_0 + h\sum_{j=0}^{s-1}\rho_{0j}(\sigma) \gamma_j(\sigma) \\ \label{sigmah}
 &\equiv& y_0 + h\sum_{j=0}^{s-1} \int_0^1 P_j(\xi)B(\sigma(\xi h))\dd \xi\int_0^1 P_j(\xi) \nabla H(\sigma(\xi h))\dd \xi,
\end{eqnarray}
in place of (\ref{infiyh}). The following results hold true.

\begin{lem}\label{grhj} With reference to (\ref{gammaro}), for any suitably regular path  $\sigma:[0,h]\rightarrow\RR^m$  one has: 
\begin{equation}\label{gammaj}
\gamma_j(\sigma)=O(h^j), \qquad \rho_{ij}(\sigma)= \rho_{ji}(\sigma) = -\rho_{ij}(\sigma)^\top=O(h^{|i-j|}),\qquad  i,j=0,1,\dots.
\end{equation}
\end{lem}
\proof See \cite[Corollary~1 and Lemma~2]{ABI2022}.\,\QED\bigskip

\begin{theo}\label{Hcons} $H(y_1)=H(y_0)$, \quad $y_1-y(h) = O(h^{2s+1})$.
\end{theo}
\proof See \cite[Theorems~1 and 2]{ABI2022}.\,\QED\bigskip

\subsection{Enforcing (\ref{teqs})}

As is clear, when applying (\ref{sigma1})--(\ref{sigmah}), to problem (\ref{poisson})-(\ref{By}), in order for $x^*$ to be reached at $t=\bar H$,
it is mandatory that the equality (\ref{teqs}) holds at the end of each integration step.  
Conversely, at $\bar H$, one would have $g(x^*)\ne0$.\footnote{Actually, $x^*$ is the order $2s$ approximation provided by the method.}
Consequently, one must have, by setting $e_m\in\RR^m$ the last unit vector,
$$
e_{m}^\top y_1 = h\,e_m^\top\int_0^1  \dot\sigma(ch)\dd c = h,
$$
i.e.,
\begin{equation}\label{teqh}
e_m^\top \int_0^1  \dot\sigma(ch)\dd c = 1.
\end{equation}
For this purpose, we specialize, for the problem at hand, the strategy used in \cite{ABI2022} for enforcing the conservation of Casimirs, thus resulting into a specific version of {\em Enhanced PHBVMs (EPHBVMs)}. Let us then consider, for a generic skew-symmetric matrix 
\begin{equation}\label{tildeB}
\tilde{B}^\top = -\tilde{B}\in\RR^{m\times m},
\end{equation}
the following modified polynomial in place of (\ref{sigma1}):\footnote{Here, we take into account that $P_0(c)\equiv1$.}
\begin{equation}\label{sigma1alfa}
\dot \sigma_\alpha(ch) =   \sum_{i,j=0}^{s-1} P_i(c)\rho_{ij}(\sigma_\alpha) \gamma_j(\sigma_\alpha) - \alpha \tilde{B}\gamma_0(\sigma_\alpha), \qquad c\in[0,1], \qquad \sigma_\alpha(0)=y_0, 
\end{equation}
with $\alpha$ a scalar to be determined. The following result holds true.

\begin{theo}\label{allalfa}
Setting 
\begin{equation}\label{y1alfa}
y_1~:=~\sigma_\alpha(h) ~\equiv~ y_0+h\sum_{j=0}^{s-1} \rho_{0j}(\sigma_\alpha) \gamma_j(\sigma_\alpha) - \alpha h \tilde{B}\gamma_0(\sigma_\alpha),
\end{equation}
one has $H(y_1)=H(y_0)$, whichever the value of $\alpha$ considered in (\ref{sigma1alfa}).
\end{theo}
\proof In fact, one has:
\begin{eqnarray*}\lefteqn{
H(y_1)-H(y_0) ~=~ H(\sigma_\alpha(h))-H(\sigma_\alpha(0)) ~=~h\int_0^1 \nabla H(\sigma_\alpha(ch))^\top\dot\sigma_\alpha(ch)\dd c}\\
&=& h\sum_{i,j=0}^{s-1} \underbrace{\int_0^1\nabla P_i(c)H(\sigma_\alpha(ch))^\top\dd c}_{=\,\gamma_i(\sigma_\alpha)^\top} \rho_{ij}(\sigma_\alpha)\gamma_j(\sigma_\alpha)
-\alpha h\underbrace{\int_0^1\nabla H(\sigma_\alpha(ch))^\top\dd c}_{=\,\gamma_0(\sigma_\alpha)^\top}\tilde{B}\gamma_0(\sigma_\alpha)\\
&=& h\sum_{i,j=0}^{s-1}\gamma_i(\sigma_\alpha)^\top\rho_{ij}(\sigma_\alpha)\gamma_j(\sigma_\alpha) ~-~\alpha h \gamma_0(\sigma_\alpha)^\top\tilde{B}\gamma_0(\sigma_\alpha)~=~0,
\end{eqnarray*}
by virtue of Lemma~\ref{grhj} and (\ref{tildeB}).\,\QED\bigskip

At this point, in order to enforce (\ref{teqs}), according to (\ref{teqh}) we require:
$$
1 ~=~\int_0^1 e_m^\top\dot\sigma_\alpha(ch)\dd c ~=~ \sum_{j=0}^{s-1} e_m^\top \rho_{0j}(\sigma_\alpha) \gamma_j(\sigma_\alpha) - \alpha  e_m^\top \tilde{B}\gamma_0(\sigma_\alpha),
$$
i.e.,
\begin{equation}\label{alfa}
\alpha = \frac{\sum_{j=0}^{s-1} e_m^\top \rho_{0j}(\sigma_\alpha) \gamma_j(\sigma_\alpha)-1}{e_m^\top \tilde{B}\gamma_0(\sigma_\alpha)}.
\end{equation}
Since $\sigma_0\equiv \sigma$, from Theorem~\ref{Hcons}, we now that the numerator is $O(h^{2s})+O(\alpha)$. Consequently, from (\ref{y1alfa}) it follows that the order of the method remains $2s$, provided that the denominator in (\ref{alfa})
is bounded away from 0. For this purpose, setting (according to (\ref{lope_y}))
$$\sigma_\alpha(ch) =: \pmatrix{c} x_\alpha(ch) \\ \omega_\alpha(ch)\endpmatrix, \quad x_\alpha(ch)\approx x(ch)\in\RR^n, \quad \omega_\alpha(ch) = \omega(ch)\in\RR, \qquad c\in[0,1],$$
and recalling that
\begin{equation}\label{gamma0}
\gamma_0(\sigma_\alpha) = \int_0^1 \nabla H(\sigma_\alpha(ch))\dd c = \pmatrix{c} \int_0^1 \nabla g(x_\alpha(ch))\dd c \\[1mm] -1\endpmatrix =: \pmatrix{c} g_0 \\[1mm] -1\endpmatrix,
\end{equation}
the choice
\begin{equation}\label{tildeB1}
\tilde B \,=\, \pmatrix{cc} O & -d_0\\ d_0^\top &0\endpmatrix, \qquad d_0 \,:=\, \frac{g_0}{ \|g_0\|_2^2},
\end{equation} 
provides 
$$e_m^\top\tilde B\gamma_0 ~=~ d_0^\top g_0~=~1.$$
In so doing, the approximation (\ref{y1alfa}) becomes, by virtue of (\ref{tildeB1}),
\begin{equation}\label{y1alfa1}
y_1~=~ y_0+h\sum_{j=0}^{s-1} \rho_{0j}(\sigma_\alpha) \gamma_j(\sigma_\alpha) - \alpha h \pmatrix{c} d_0 \\[1mm] 1\endpmatrix,
\end{equation}
with 
\begin{equation}\label{alfa1}
\alpha = \sum_{j=0}^{s-1} e_m^\top \rho_{0j}(\sigma_\alpha) \gamma_j(\sigma_\alpha)-1. 
\end{equation}

\subsection{Discretization}

As is clear, the Fourier coefficients
$$
\gamma_j(\sigma_\alpha)     = \int_0^1 P_j(\xi)\nabla H(\sigma_\alpha(\xi h))\dd\xi, \quad
\rho_{ij}(\sigma_\alpha)         = \int_0^1 P_i(\xi)P_j(\xi)B(\sigma_\alpha(\xi h))\dd\xi, \quad i,j=0,\dots,s-1,
$$
need to be numerically computed. For this purpose, we use a Gauss-Legendre formula of order $2k$, with abscissae and weights
$(c_i,b_i)$, $i=1,\dots,k$. In so doing, we obtains a new polynomial approximation, say $u_\alpha$, in place of $\sigma_\alpha$,
\begin{equation}\label{u1alfa}
\dot u_\alpha(ch) =   \sum_{i,j=0}^{s-1} P_i(c)\hat\rho_{ij}(u_\alpha) \hat\gamma_j(u_\alpha) - \alpha h \pmatrix{c} \hat d_0 \\[1mm] 1\endpmatrix, \qquad c\in[0,1],
\end{equation}
where, setting as before,
$$u_\alpha(ch) =: \pmatrix{c} x_\alpha(ch) \\ \omega_\alpha(ch)\endpmatrix, \quad x_\alpha(ch)\approx x(ch)\in\RR^n, \quad \omega_\alpha(ch) = \omega(ch)\in\RR, \qquad c\in[0,1],$$
$\hat d_0$ is defined (compare with (\ref{tildeB1})) as
\begin{equation}\label{hg0} 
\hat d_0 = \frac{\hat g_0}{\|\hat g_0\|_2^2}, \qquad \hat g_0 = \sum_{\ell=1}^k b_\ell\nabla g(u_\alpha(c_\ell h)),
\end{equation}
and we use the (generally) approximate Fourier coefficients
\begin{eqnarray}\label{hfourier}
\hat\gamma_j(u_\alpha)     &=& \sum_{\ell=1}^k b_\ell P_j(c_\ell)\nabla H(u_\alpha(c_\ell h)), \\ \nonumber
\hat\rho_{ij}(u_\alpha)         &=& \sum_{\ell=1}^k b_\ell P_i(c_\ell)P_j(c_\ell)B(u_\alpha(c_\ell h)), \qquad i,j=0,\dots,s-1.
\end{eqnarray}
At last, $\alpha$ is defined as (compare with (\ref{alfa1})):
\begin{equation}\label{alfa1d}
\alpha = \sum_{j=0}^{s-1} e_m^\top \hat\rho_{0j}(u_\alpha) \hat\gamma_j(u_\alpha)-1. 
\end{equation} 
Setting, as usual (compare with (\ref{y1alfa})),
\begin{equation}\label{y1u}
y_1 ~:=~y_0 +h\sum_{j=0}^{s-1} \hat\rho_{0j}(u_\alpha) \gamma_j(\sigma_\alpha) - \alpha h \pmatrix{c} \hat d_0 \\[1mm] 1\endpmatrix,
\end{equation}
the following results follow.

\begin{theo}\label{kges} $\forall k\ge s\,:\, y_1-y(h) = O(h^{2s+1}).$\end{theo}
\proof See \cite[Theorem~10]{ABI2022}.\,\QED\bigskip

\begin{theo}\label{gpol} With reference to (\ref{hg0})--(\ref{y1u}), if  $g\in\Pi_\nu$ and  $\nu\le\frac{2k}s$, then
\begin{equation}\label{sono_esatti}
\hat g_0 = g_0 \equiv \int_0^1 \nabla g(x_\alpha(ch))\dd c, \quad \hat\gamma_j(u_\alpha) = \gamma_j(u_\alpha) \equiv
\int_0^1 P_j(\xi)\nabla H(u_\alpha(\xi h))\xi, \quad j=0,\dots,s-1.
\end{equation}
Consequently, $H(y_1) = H(y_0)$.
\end{theo}
\proof The first statement follows from the fact that the integrands in (\ref{sono_esatti}) are polynomials of degree at most $(\nu-1)s+s-1 = \nu s-1\le 2k-1$. Energy conservation is then proved with similar steps as in the proof of Theorem~\ref{allalfa}, by formally replacing $\sigma_\alpha$ with $u_\alpha$.\QED\bigskip

\begin{rem}\label{gnotpol} In the case where $g$ is not a polynomial, or is a polynomial but the hypotheses of the previous Theorem~\ref{gpol} are not fulfilled,
form \cite[Theorem~9]{ABI2022} it follows that $$H(y_1)=H(y_0)+O(h^{2k+1}).$$
Consequently, a {\em practical} energy-conservation can always be gained, provided that $k$ is chosen large enough so that the energy error falls below the round-off error level.
\end{rem}

\begin{defi} According to \cite[Definition~2]{ABI2022}, we shall refer to the numerical method defined by (\ref{hg0})--(\ref{y1u}) as  the EPHBVM$(k,s)$ method.
\end{defi}

\begin{rem}\label{keqs} We observe that, when $k=s$, the EPHBVM$(s,s)$ method naturally satisfies the constraint (\ref{teqs}). Consequently,  $\alpha=0$ and the method coincides with the symplectic $s$-stage Gauss-Legendre collocation method used for solving (\ref{poisson})-(\ref{By}).\end{rem}

For sake of completeness, let us sketch the vector form of the EPHBVM$(k,s)$ method, which can be derived by slightly adapting the arguments in \cite[Section~4.1]{ABI2022}.
For this purpose, let us define the matrices (see (\ref{orto}))
\begin{eqnarray*}
&&\P_s = \pmatrix{c} P_{j-1}(c_i)\endpmatrix_{\scriptsize\begin{array}{l}i=1,\dots,k\\ j=1,\dots,s\end{array}}, ~ \I_s = \pmatrix{c} \int_0^{c_i}P_{j-1}(\xi)\dd\xi\endpmatrix_{\scriptsize\begin{array}{l}i=1,\dots,k\\ j=1,\dots,s\end{array}} \in\RR^{k\times s}, \\
&&\Omega = \diag(b_1,\dots,b_k),
\end{eqnarray*}
with $(c_i,b_i)$, $i=1,\dots,k$, the abscissae and weights of the Gauss-Legendre quadrature, and the vectors (see (\ref{hfourier}))
$$\bfe = \pmatrix{c} 1\\ \vdots \\ 1\endpmatrix\in\RR^k, \qquad \bfphi = \pmatrix{c} \phi_0\\ \vdots \\ \phi_{s-1}\endpmatrix, \quad \phi_i = \sum_{j=0}^{s-1} \hat\rho_{ij}(u_\alpha)\hat\gamma_j(u_\alpha)\in\RR^m,\quad i=0,\dots,s-1.$$
 We also set, being\,\footnote{According to (\ref{lope_y}), \,$x_i\in\RR^n$\, and\, $\omega_i\in\RR$.}
 \begin{equation}\label{Yixi}
 Y \equiv \pmatrix{c} Y_1\\ \vdots \\ Y_k\endpmatrix\in\RR^{km},\qquad Y_i\equiv \pmatrix{c} x_i \\ \omega_i\endpmatrix,  \quad i=1,\dots,k,
 \end{equation}
 the stages of the method, and (see (\ref{poisson})-(\ref{By})),
 $$\B(Y) = \pmatrix{ccc} B(Y_1)\\ &\ddots\\ &&B(Y_k)\endpmatrix\in\RR^{km\times km},$$
 one then obtains the discrete problem
 \begin{equation}\label{dispro}
 \F(\bfphi,\alpha) := \pmatrix{c}
 \bfphi-\P_s^\top\Omega\otimes I_m \B(Y) \cdot \P_s\P_s^\top\Omega\otimes I_m\nabla H(Y)\\[1mm]
 \alpha - e_m^\top\phi_0+1\endpmatrix = \bfzero\in\RR^{sm+1},
 \end{equation}
 where, with reference to (\ref{Yixi}), and setting ~$\bfc=(c_1,\dots,c_k)^\top$~ the vector of the abscissae:
\begin{eqnarray*}
Y &=& \bfe\otimes y_0 + h\I_s\otimes I_m\bfphi - \alpha h\bfc\otimes \pmatrix{c}\hat d_0\\1 \endpmatrix,\\
\hat d_0 &=& \frac{\hat g_0}{\|\hat g_0\|_2^2},\\
\hat g_0 &=& \sum_{i=1}^k b_i\nabla g(x_i).
\end{eqnarray*}

\begin{rem} We observe that the above discrete problem (\ref{dispro}) has, remarkably, (block) dimension $s$, independently of the considered value of $k$ \cite{ABI2022}. Moreover, it induces a straightforward fixed-point iteration, which converges for all sufficiently small timesteps $h$, under regularity assumptions on $f$ and $g$. This iteration will be used for the numerical tests, even though Newton-type procedures, obtained adapting those defined in \cite{BFCI2014,BIT2011} for HBVMs (see also \cite{BIM2015,BM2002,BM2009}), could be also considered.
\end{rem}

\section{Numerical tests}\label{numtest}

In this section we present a few numerical tests, concerning the solution of  one-sided event location problems, aimed at assessing the theoretical findings. For each problem we prescribe the function $f(x)$ in (\ref{lope_ode}) and the event function $g(x)$ (\ref{lope_h}), along with the starting point of the trajectory. All the numerical tests have been implemented in Matlab (R2020a) on a 3 GHz Intel Xeon W 10 core computer with 64 GB of memory.

\paragraph{Example~1} 
The first test problem, taken from \cite[Example 5.1 (a)]{DL2015}, is defined by
\begin{equation}\label{ex1_f}
f(x) = \pmatrix{c} x_2\\[1mm] \frac{1}{1.2-x_2}-x_1\endpmatrix, 
\end{equation}
and by the event function  
\begin{equation}\label{ex1_g}
g(x) = x_1+x_2-0.4.
\end{equation}
We choose the initial point 
\begin{equation}\label{ex1_0}
x(0) = (-0.2,\,-0.2)^\top, 
\end{equation}
providing the value $\bar H\equiv -g(x(0)) = 0.8$ in (\ref{lope_H}). Since the event function is linear, any EPHBVM$(s,s)$ method (i.e., the $s$-stage Gauss collocation method) has to provide an order $2s$ approximation to the event point belonging to the event set $\Sigma$  (\ref{Sset}).This is confirmed by the numerical tests listed in Table~\ref{ex1_tab}, obtained by using the timesteps
\begin{equation}\label{ex1_h}
h_n := \frac{\bar H}{10\cdot 2^n}
\end{equation}
for solving the associated Poisson problem (\ref{poisson})-(\ref{By}). 
In the table, we have denoted by $x_n^*$ the approximation to the event point obtained with the timestep $h_n$, and with $e_n^*$ the corresponding error (numerically estimated). As one may see, all approximations belong to the event set $\Sigma$ (since $g(x_n^*)$ is of the order of the round-off error level) and converge to the event point with the correct order (the last two approximations for $s=3$ practically coincide).

\begin{table}[t]
\caption{numerical results for problem (\ref{ex1_f})--(\ref{ex1_0}) solved by the EPHBVM$(s,s)$ method with timestep (\ref{ex1_h}).}
\label{ex1_tab}\centering
\begin{tabular}{|c|rcc|rcc|rcc|}
\hline
&\multicolumn{3}{c|}{$s=1$} &\multicolumn{3}{c|}{$s=2$} &\multicolumn{3}{c|}{$s=3$}\\
\hline
$n$ & $g(x_n^*)$ & $e_n^*$ & rate & $g(x_n^*)$ & $e_n^*$ & rate & $g(x_n^*)$ & $e_n^*$ & rate\\
\hline 
  0 & -1.11e-16 & --- & --- & -1.67e-16 &--- & --- & -1.05e-15 &--- & --- \\ 
  1 & -5.55e-17 &1.99e-04 & --- & -1.67e-16 &1.43e-07 & --- & -9.44e-16 &6.09e-10 & --- \\ 
  2 & 1.11e-16 &4.98e-05 & 2.0 & 0.00e+00 &9.33e-09 & 3.9 & -9.99e-16 &1.06e-11 & 5.8 \\ 
  3 & -1.67e-16 &1.25e-05 & 2.0 & -2.22e-16 &5.90e-10 & 4.0 & -8.88e-16 &1.72e-13 & 6.0 \\ 
  4 & -5.55e-17 &3.11e-06 & 2.0 & -1.67e-16 &3.70e-11 & 4.0 & -1.22e-15 &2.60e-15 & 6.0 \\ 
  5 & 2.22e-16 &7.78e-07 & 2.0 & -1.11e-16 &2.31e-12 & 4.0 & -1.11e-15 &1.39e-16 & *** \\ 
\hline
\end{tabular}
\end{table}

\paragraph{Example~2}  
Next, we consider the problem defined by (\ref{ex1_f}), with a nonlinear (though smooth) event function
\begin{equation}\label{ex2_g}
g(x) = 20x_1+x_2-20\sin x_1 - 0.4,
\end{equation}
and initial point
\begin{equation}\label{ex2_0}
x(0) = (0,\,-0.2)^\top, 
\end{equation}
providing a value $\bar H=0.6$. If we solve the associated Poisson problem (\ref{poisson})-(\ref{By}) by using the EPHBVM$(s,s)$ and EPHBVM$(4,s)$ methods, $s=1,2,3$, with a timestep $h=\bar H/10$, we see that, though $g(x)$ is non-polynomial, the latter methods correctly reaches the event set $\Sigma$, as one infers from the results listed in Table~\ref{ex2_tab}.

\begin{table}[t]
\caption{numerical results for problem (\ref{ex1_f}) and (\ref{ex2_g})-(\ref{ex2_0}) solved by the EPHBVM$(k,s)$ method with timestep $h=\bar H/10$.}
\label{ex2_tab}\centering
\begin{tabular}{|c|c|r|r|}
\hline
              & $s$  & EPHBVM$(s,s)$ & EPHBVM$(4,s)$\\
\hline 
              & 1 &  1.1148e-05 & 1.1102e-16 \\ 
$g(x^*)$ & 2 & -1.4687e-08 & 1.1102e-16 \\ 
              & 3 &  -7.8148e-11& -2.2204e-16 \\ 
\hline
\end{tabular}
\end{table}
   
\begin{figure}[t]
\centering
\includegraphics[width=12cm]{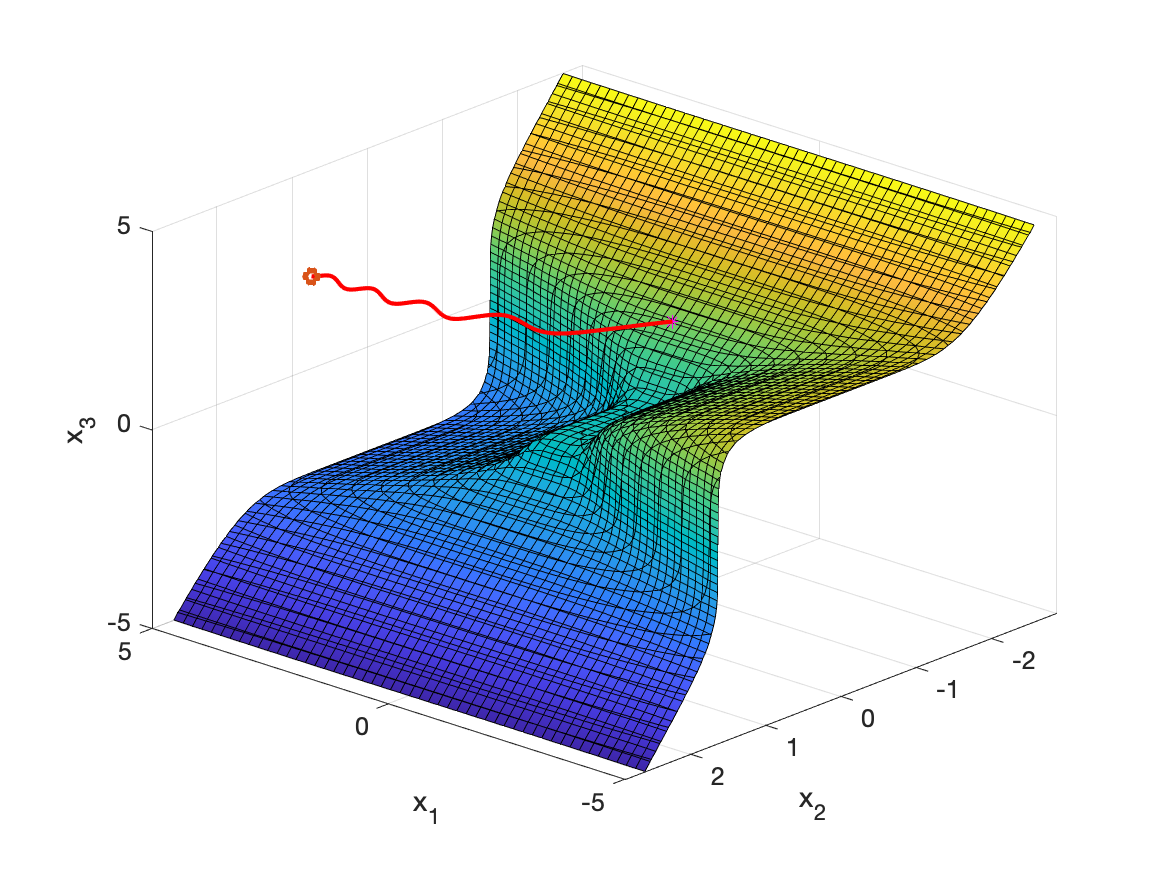}
\caption{event set (\ref{Sset})  for problem (\ref{ex3_f})-(\ref{ex3_g}), along with the trajectory reaching it computed by using the EPHBVM(11,3) method.}
\label{lopez3_3}
\end{figure}

\paragraph{Example~3}  
At last, let us consider the problem defined by:
\begin{equation}\label{ex3_f}
f(x) = \pmatrix{c} \frac{1}{1.2+\sin x_2} \\[1.5mm] \frac{1}{1.2-\cos x_1}\\[1.5mm]   1+\cos \|x\|_2^2\endpmatrix,  \qquad x(0) = \pmatrix{c} 3\\ 2\\ 4\endpmatrix,
\end{equation}
with the polynomial event function\,\footnote{The scaling factor $10^{-3}$ in (\ref{ex3_g}) is introduced to have a more compact graphical representation.}
\begin{equation}\label{ex3_g}
g(x) = 10^{-3}\left( x_1^3+4x_2^7+x_3^5\right).
\end{equation}
For the chosen initial point, one obtains $\bar H = 1.563$. For this problem, any EPHBVM$(k,s)$ method, with
$$k\ge\frac{7s}2,$$
turns out to be energy-conserving for the associated Poisson problem (\ref{poisson})-(\ref{By}), and therefore, at $t=\bar H$ the trajectory exactly reaches the event point lying on the event set $\Sigma$. We use the EPHBVM(11,3) method with timestep $h=10^{-3}\bar H$, thus reaching the (approximation of the) event point $x^*$ for which $g(x^*)\approx -5.8\cdot 10^{-16}$. The event set $\Sigma$, along with the computed trajectory, are depicted in Figure~\ref{lopez3_3}. 
 For comparison, the EPHBVM(3,3) (i.e., the 3-stage Gauss-Legendre method), using the same timestep, reaches a point $\tilde x$ for which $g(\tilde x)\approx 4.7\cdot 10^{-8}$ and, moreover, $\|\tilde x-x^*\|_2\approx 4.8\cdot 10^{-6}$.

\section{Conclusions}\label{fine}

In this paper, starting from the methodology introduced in \cite{DL2015}, we have introduced a direct method for numerically solving the problem of one-sided event location. The proposed approach is based on a suitable modification of recently derived energy-conserving methods for Poisson problems \cite{ABI2022}, specifically tailored for the problem at hand. The methods exactly reach the event set, in the case where the event function is a polynomial. Actually, they can be effectively used also in the non-polynomial case, provided that the event function is regular enough. Numerical examples confirm the theoretical findings.

\subsection*{Acknowledgements} The authors wish to thanks the {\em mrSIR project} \cite{mrsir} for the financial support.

\subsection*{Declaration of interest} The authors declare no competing interest.

\end{document}